\documentclass{amsart}

\usepackage{graphicx}

\newtheorem{theorem}{Theorem}[section]
\newtheorem{lemma}[theorem]{Lemma}
\newtheorem{corollary}[theorem]{Corollary}

\theoremstyle{definition}
\newtheorem{definition}[theorem]{Definition}

\theoremstyle{remark}
\newtheorem{remark}[theorem]{Remark}

\numberwithin{equation}{section}

\def\qed{\hfill $\square$}

\newcommand{\eps}{\varepsilon}

\begin{document}

\title[Stable manifolds for periodically perturbed maps]{
Stable manifolds for periodically perturbed maps
}


\author{Matthew Williams
}
\address{Department of Mathematical Sciences, University of Texas at Dallas, 75080, Richardson, USA
}
\curraddr{}
\email{matthew.williams3@utdallas.edu
}
\thanks{}

\author{Oleg Makarenkov
}
\curraddr{}
\email{
makarenkov@utdallas.edu
}
\thanks{}


\date{}

\dedicatory{}

\commby{}

\begin{abstract} We prove that if a certain entry in the map of the Hadamard-Perron theorem is $T$-periodic in one of the variables, then the stable manifold guaranteed by the Hadamard-Perron theorem is a graph of a $T$-periodic function. As an application, we extend the classical Levinson's result about the occurrence of an attracting closed invariant curve near a stable cycle of a system of autonomous equations under periodic perturbations to hybrid differential equations.  
\end{abstract}

\maketitle


\bibliographystyle{amsplain}

\section{Introduction} 

\begin{figure}
\includegraphics[scale=0.35]{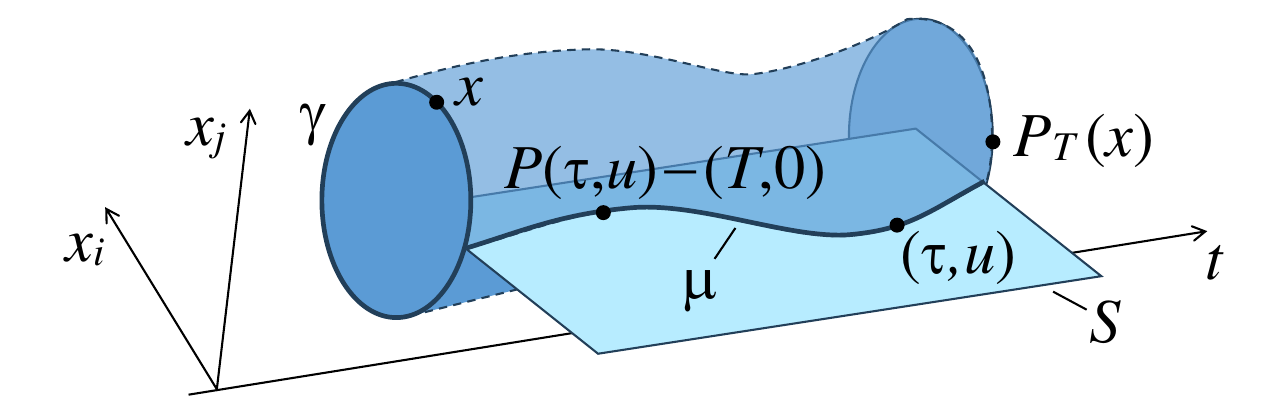} 
\label{fig1} 
\caption{The two types of attracting invariant manifolds $\gamma$ and $\mu$ (of Poincar\'e maps $P_T$ and $P$) given by Levinson's theorem. The left and right discs correspond to $t=0$ and $t=T$.  
}
\end{figure}

\noindent 
As showed by Levinson \cite{Levinson}, the dynamics of $T$-periodically perturbed self-oscillating systems is one-dimensional (given by a continuous map on a circle). Let $P_T$ be the stroboscopic map that, for a solution $x(t)$ of the perturbed system, maps every initial condition $x(0)$ to the point $x(T)$. To establish the existence of a one-dimensional invariant curve $\gamma$ for $P_T$,  the work \cite{Levinson} applies the Cacciopoli fixed point theorem and obtains $\gamma$ as a fixed point in a functional space of curves. A visualization of Levinson's result is given in Fig.~\ref{fig1}, where solutions originating at $\{0\}\times \gamma$ land at $\{T\}\times \gamma$ after time $T$ forming a $T$-periodic cylinder. Another interpretation of Levinson's result can be obtained by drawing a cross-section $S$ 
through $\{0\}\times \gamma(\tau)$ and $\{T\}\times \gamma(\tau)$ transversely to the surface of the cylinder as shown at Fig.~\ref{fig1}. The result \cite{Levinson} implies that the Poincar\'e map $P$ induced by $S$ admits a $T$-periodic  attracting invariant curve $\mu$ as well. We discover that such a $\mu$ exists regardless of whether the differential equations are smooth or nonsmooth along $S$, thus allowing $S$ to be a switching surface, which is the case e.g. in robotic locomotion \cite{Makarenkov,Westervelt}.

\vskip0.2cm

\noindent Stable invariant manifolds for nonsmooth differential equations have been shown to exist e.g. in the case of a perturbed equilibrium lying on the switching surface \cite{Weiss}. When the reduced system admits more complex attractors, such as e.g. a cycle,  the available results about non-standard behavior of invariant manifolds concerned the case where the switching surface is not transversal to the cycle \cite{Kryzhevich,Szalai}. However, whether or not the case of transversal intersection resembles the smooth dynamics (given by Levinson's theorem) remained an open problem.

\vskip0.2cm

\noindent The present paper casts the problem of the existence of a closed invariant curve $\mu$ as a problem of periodicity in the invariant manifolds coming from Hadamard-Perron theorem. Indeed, when $\eps=0$ the cylinder of Fig.~\ref{fig1} consists of the cycle of the reduced autonomous system translated along the $t$-axis from $t=0$ to $t=T$. In other words, when $\eps=0$, the invariant curve $\mu$ is a straight line $L$ obtained from the fixed point $u_0$ of $P(0,u)$ (corresponding to the cycle of the autonomous system) translated from $t=0$ to $t=T$. Assuming that fixed point $u_0$ is stable, the occurrence of an invariant curve $\mu$ of $P$ near $L$ comes from Hadamard-Perron theorem after appropriate technical constructions, see the summary in Remark~\ref{C1C3}. The main contribution of the present paper is establishing of an appropriate periodicity in the invariant manifolds in Hadamard-Perron theorem (that gives us periodicity of $\mu$, thus one dimensionality of the dynamics).

\vskip0.2cm

\noindent In the field of hybrid differential equations, the result of the present paper can be viewed as the proof of the existence of an aperiodic gait \cite{Yang} in robotic walking under time-periodic disturbances
 (an open problem listed in  \cite{Grizzle}). 

\vskip0.2cm

\noindent The paper is organized as follows. In  Section~\ref{nextsec} we derive an appropriate modification of the Hadamard-Perron theorem. We recall the reader that assuming that the map $f:\mathbb{R}^{k_1}\times\mathbb{R}^{k_2}\to\mathbb{R}^{k_1+k_2}$ is a global diffeomorphism of the form
\begin{equation}\label{f-original}
  f(x,y)=(Ax+\alpha(x,y),By+\beta(x,y)),
\end{equation}
with $\|A^{-1}\|\le\mu^{-1}$, $\|B\|\le\lambda$ for some $\lambda<\mu,$ and
$
   \alpha(0)=\beta(0)=0,
$
one of the conclusions of the  Hadamard-Perron theorem (Katok-Hasselblatt \cite[Theorem~6.2.8]{Katok}) is the existence  of a curve $\varphi:\mathbb{R}^{k_1}\to\mathbb{R}^{k_2}$ which is invariant in the sense that
$$   f(\mathrm{graph}(\varphi))=\mathrm{graph}(\varphi).
$$
Because the Poincar\'e map $P(\tau,u)$ is seldom a global homeomorphism,
Section~\ref{nextsec}  replaces the assumption for $f$ to be a global diffeomorphism by boundedness of $x\mapsto\alpha(x,y)$ and $x\mapsto\beta(x,y)$  and contractivity of $y\mapsto\beta(x,y)$, which allows us to establish (Theorem~\ref{thm1}) the existence, attractivity, and boundedness of $\varphi$. The periodicity of $\varphi$ is then established in Section~\ref{Additional} (Theorem~\ref{corperiodicity}) by strengthening the boundedness of $\alpha$ and $ \beta$ to   periodicity.  An application to periodic perturbations of hybrid differential equations with a stable limit cycle is given in Section~\ref{sec-application}.  

\section{A version of the Hadamard-Perron theorem without the global diffeomorphism assumption}\label{nextsec}

\noindent With  applications to Poincar\'e maps $P(\tau,u)$ of differential equations in mind we focus on  $A=I$, i.e. consider \begin{equation}\label{Pformula0}
f_{\omega,\eps}(x,y)=\left(x+\omega \alpha(\omega,\eps,x,y),\beta(\omega,\eps,x,y)\right),
\end{equation}
 where $\beta(\omega,\eps,x,y)$ includes $By$ and $\beta(x,y)$ of (\ref{f-original}) and parameter $\omega$ is introduced for technical purposes.
In what follows, $\mathbb{D}^k$ denotes the closed unit disc in $\mathbb{R}^k$, i.e. 
$$
\mathbb{D}^k=\left\{x\in\mathbb{R}^k:\|x\|\le 1\right\}.
$$
Accordingly, $r\mathbb{D}^k$ is the closed disc of radius $r.$ 

\vskip0.2cm

\noindent To shorten the text, we will often identify vectors rows and vectors columns. For example, the notation $(\omega,\eps,x,y)$ will often mean $(\omega,\eps,x,y)^T$ when it doesn't cause a confusion (we will add the ``$T$" when confusion arises).

\begin{definition} Assume that $\alpha\in C^2(\mathbb{R}\times\mathbb{R}\times\mathbb{R}^{k_1}\times r_1\mathbb{D}^{k_2},\mathbb{R}^{k_1}),$ $\beta\in C^2(\mathbb{R}\times\mathbb{R}\times\mathbb{R}^{k_1}\times r_1\mathbb{D}^{k_2},\mathbb{R}^{k_2})$ for some $r_1>0$. If a sequence $\{(x_j,y_j)\}_{j=0}^\infty\subset\mathbb{R}\times r_1\mathbb{D}^{k_2}$ satisfies  
$(x_{j},y_{j})=f_{\omega,\eps}(x_{j-1},y_{j-1}),$ $j\in\mathbb{N},$ then the trajectory of $f_{\omega,\eps}$ with the initial the condition $(x_0,y_0)$ is defined and is the sequence $\{(x_j,y_j)\}_{j=0}^\infty.$
\end{definition}


\begin{theorem}\label{thm1} Assume that

\begin{enumerate}
\renewcommand{\theenumi}{A\arabic{enumi}} \renewcommand{\labelenumi}{\theenumi)}
\item \label{A-regularity} $\alpha\in C^2(\mathbb{R}\times\mathbb{R}\times\mathbb{R}^{k_1}\times r_1\mathbb{D}^{k_2},\mathbb{R}^{k_1}),$ $\beta\in C^2(\mathbb{R}\times\mathbb{R}\times\mathbb{R}^{k_1}\times r_1\mathbb{D}^{k_2},\mathbb{R}^{k_2}),$ for some $r_1>0$, and when 
$\omega,\eps,y$ take values from a bounded set, 
the functions 
$\alpha(\omega,\eps,x,y)$ 
and 
$\beta(\omega,\eps,x,y)$ and all their partial derivatives to order 2 
are bounded in $x\in\mathbb{R}^{k_1}.$
\item \label{A-zero} $\alpha(\omega,0,x,y)$ and $\beta(\omega,0,x,y)$ are independent of $\omega$ and $x$, and $\beta(0)=0$;
\item \label{A-stability} the matrix $\beta_y(0)$ is invertible and there exists $q\in[0,1)$ such that $\|\beta(0,0,0,y_1)-\beta(0,0,0,y_2)\|\le q\|y_1-y_2\|,$ $y_1,y_2\in r_1\mathbb{D}^{k_2}.$ 
\end{enumerate}
Then there exist $\eps_0>0$ and  $0<r_0<r_1$, such that for any $\omega\in[0,1]$ and any $\eps\in(-\eps_0,\eps_0)$ there exists a unique $\varphi_{\omega,\eps}\in C^1(\mathbb{R}^{k_1},r_0\mathbb{D}^{k_2})$ satisfying
\begin{enumerate}
\renewcommand{\theenumi}{C\arabic{enumi}}
\renewcommand{\labelenumi}{\theenumi)}

    \item\label{C1}
$
f_{\omega,\eps}(\mathrm{graph}(\varphi_{\omega,\eps}))=\mathrm{graph}(\varphi_{\omega,\eps});
$
\item \label{C2} the trajectory $\{(x_j,y_j)\}_{j=0}^\infty$ of  dynamical system $f_{\omega,\eps}$ with initial condition $(x_0,y_0)\in\mathbb{R}^{k_1}\times r_0\mathbb{D}^{k_2}$ is defined for all $j\in\mathbb{N}$ and 
$$
   \|y_j-\varphi_{\omega,\eps}(x_j)\|\to 0\quad {\rm as}\quad j\to\infty,
$$
 uniformly with respect to the initial condition $(x_0,y_0)$.
\end{enumerate}
Furthermore
\begin{enumerate}
\renewcommand{\theenumi}{C\arabic{enumi}}
\renewcommand{\labelenumi}{\theenumi)}
\setcounter{enumi}{2}
\item \label{C3} the function
$\varphi_{\omega,\eps}$ depends continuously on  $\omega\in[0,1]$ and $\eps\in(-\eps_0,\eps_0)$, and $\varphi_{\omega,0}(x)\equiv 0.$
In particular, $x\mapsto \varphi_{\omega,\eps}(x)$ is bounded on $\mathbb{R}^{k_1}$. 
\end{enumerate}

\end{theorem}


\noindent To prove the theorem we embed map $f_{\omega,\eps}$ into the map
\begin{equation}\label{Flambda}
  F_{\lambda}(\omega,\eps,x,y)=\left(A_\lambda(\omega,\eps,x)^T+\tilde \alpha_{\lambda}(\omega,\eps,x,y),B y+\tilde \beta_{\lambda}(\omega,\eps,x,y)\right),
\end{equation}
where
$$
\begin{array}{ll}
  A_\lambda
  =\left(\begin{array}{ccc}1 & 0 & 0 \\ 0 & 1 & 0\\  \lambda\alpha(0) & 0 & I\end{array}\right),&  \tilde\alpha_{\lambda}(\omega,\eps,x,y)=\left(\begin{array}{c} 0 \\ 0 \\ \lambda\omega\left(\alpha(\lambda\omega,\lambda^2\eps,x,\lambda y)-\alpha(0)\right)\end{array}\right),\\ B=\beta_y(0),
  & \tilde\beta_{\lambda}(\omega,\eps,x,y)=\dfrac{1}{\lambda}\beta(\lambda\omega,\lambda^2\eps,x,\lambda y)-\beta_y(0)y.
\end{array}
$$

\begin{remark}\label{Ff} \rm It holds that the first and second rows of $F_{\lambda}(\omega,\eps,x,y)$ are $\omega$ and $\eps$ respectively, and  $$
F_{\lambda}(\omega,\eps,x,y)=(\omega,\eps,\bar x,\bar y),\quad\mbox{if and only if}\quad f_{\lambda\omega,\lambda^2\eps}(x,\lambda y)=(\bar x,\lambda\bar y).$$
\end{remark}

\noindent We now want to extend $F_\lambda$ to a global $C^1$ diffeomorphism through a suitable bump function. The next lemma evaluates the support  (domain) of the bump function. 


\vskip0.2cm


 \begin{lemma}\label{magnitude}  
Let $\tilde\alpha_\lambda$ and $\tilde\beta_\lambda$ be the nonlinearities from (\ref{Flambda}) and let $D\tilde\alpha_\lambda$ and $D\tilde\beta_\lambda$ be the corresponding derivatives. Under assumptions (\ref{A-regularity}) and (\ref{A-zero}) of Theorem~\ref{thm1}, for any $\delta>0$ there exists $\lambda_0>0$ such that, for any $\lambda\in(0,\lambda_0],$
\begin{eqnarray} && \hskip-0.5cm
\|\tilde\alpha_\lambda(\omega,\eps,x,y)\|\hskip-0.05cm+\hskip-0.05cm\|D\tilde\alpha_\lambda(\omega,\eps,x,y)\|\le\delta,\ 
\|\tilde\beta_\lambda(\omega,\eps,x,y)\|\hskip-0.05cm+\hskip-0.05cm\|D\tilde\beta_\lambda(\omega,\eps,x,y)\|\le\delta, \nonumber \\
&& \qquad \omega  \in\left[-\dfrac{1}{\lambda},\dfrac{2}{\lambda}\right], \ \eps\in(-r,r),\ x\in\mathbb{R}^{k_1}, \ \|y\|\le r_1.  \label{required}
 \end{eqnarray}
      \end{lemma}

\noindent {\bf Proof.} The required estimate for $\tilde\alpha_\lambda$ comes from the boundedness assumption of (\ref{A-regularity}). For the rest of this proof, $\beta_\eps(\omega,\eps,x,y)$ and $\beta_y(\omega,\eps,x,y)$ are the partial derivatives of $\beta(\omega,\eps,x,y)$ with respect to $\eps$ and $y$ respectively.   
To get the required estimate for $\tilde\beta_\lambda$ we first use (\ref{A-zero}) 
to represent 
$$ \dfrac{1}{\lambda}\beta(\lambda\omega,\lambda^2\eps,x,\lambda y)-\beta_y(0)y=\dfrac{1}{\lambda}\beta(\lambda\omega,\lambda^2\eps,x,\lambda y)-\dfrac{1}{\lambda}\beta(\lambda\omega,0,x,0)-\beta_y(\lambda\omega,0,x,0)y.
$$
Taylor Theorem with  remainder in the integral form (see e.g. \cite[Taylor Theorem of Ch.~XIII, \S6]{Lang}) yields 
\begin{equation}\label{betaestimate}
\begin{array}{l}
\dfrac{1}{\lambda}\beta(\lambda\omega,\lambda^2\eps,x,\lambda y)-\beta_y(0)y=\int\limits_0^1\beta_\eps(\lambda\omega,(t\lambda)^2\eps,x,t\lambda y)2\lambda\eps\, dt+\\ \hskip3cm+\int\limits_0^1\left(\beta_y(\lambda\omega,(t\lambda)^2\eps,x,t\lambda y)y-\beta_y(\lambda\omega,0,x,0)y\right)dt.
\end{array}
\end{equation}
The first integral in the sum can be made arbitrary small by reducing $\lambda_0>0$ due to the boundedness of $\beta_\eps$ and the second term in the sum can be made arbitrary small by $C^1$ regularity of $\beta$. The fact that  $\|D\tilde\beta_\lambda(\omega,\eps,x,y)\|$ can be made as small as necessary by further reduction of $\lambda_0>0$ follows by taking partial derivaties of (\ref{betaestimate}) and by using the boundedness and continuity of 2nd order partial derivatives of $\beta$ as insured by (\ref{A-regularity}). \qed


\vskip0.2cm

\noindent Based on Lemma~\ref{magnitude}, we introduce a 
 $ C^1(\mathbb{R}\times\mathbb{R}\times\mathbb{R}^{k_2},[0,1])$ bump function 
\begin{equation}\label{bumpfunction}
  \Psi(\omega,\eps,y)=\left\{\begin{array}{cc}
  1, & (\omega,\eps,y)\in[0,1]\times\left(-\dfrac{r_1}{2},\dfrac{r_1}{2}\right)\times\dfrac{r_1}{2}\mathbb{D}^{k_1},\\
  0, & \mathbb{R}\times\mathbb{R}\times\mathbb{R}^{k_2}\backslash[-1,2]\times\left(-r_1,r_1\right)\times r_1\mathbb{D}^{k_2},\end{array}\right.
\end{equation}
and vanish the the nonlinear terms $\tilde\alpha_\lambda$ and $\tilde\beta_\lambda$ in (\ref{Flambda}) as follows:
$$
G_\lambda(\omega,\eps,x,y)=  \left(\begin{array}{cc}A_\lambda & 0 \\ 0 & B\end{array}\right)(\omega,\eps,x,y)^T+\Psi(\lambda\omega,\eps,y)\left(\begin{array}{c}\tilde \alpha_{\lambda}(\omega,\eps,x,y) \\ \tilde \beta_{\lambda}(\omega,\eps,x,y)\end{array}\right).
$$
\begin{corollary}\label{Glambdaproof}
    Under assumptions 
    of Theorem~\ref{thm1}, there exists $\lambda_0>0,$ such that for any $\lambda\in(0,\lambda_0]$, $G_\lambda$ is a surjective $C^1$-diffeomorphism from $\mathbb{R}\times\mathbb{R}\times\mathbb{R}^{k_1}\times\mathbb{R}^{k_2}$ to $\mathbb{R}\times\mathbb{R}\times\mathbb{R}^{k_1}\times\mathbb{R}^{k_2}$.
\end{corollary}

\noindent The proof of Corollary~\ref{Glambdaproof} uses a known result saying that a small nonlinear perturbation of a linear invertible map is a diffeomorphism, see Theorem~\ref{diffthm} in  Appendix. With this result, 
 Corollary~\ref{Glambdaproof}  comes  by noticing the following two remarks.
\begin{remark} The $(k_1+k_2+2)\times (k_1+k_2+2)$-matrix $\left(\begin{array}{cc}A_\lambda & 0 \\ 0 & B\end{array}\right)$ is invertible for all $\lambda>0$ sufficiently small because it is invertible at $\lambda=0$ by (\ref{A-stability}). \end{remark}
\begin{remark}
By Lemma~\ref{magnitude}, for any $\delta>0$ there exists $\lambda_0>0$ such that for $\lambda\in(0,\lambda_0],$ the $C^1$-norms of functions 
$$
(\omega,\eps,x,y)\mapsto \Psi(\lambda\omega,\eps,y)\tilde\alpha_\lambda(\omega,\eps,x,y),\quad (\omega,\eps,x,y)\mapsto \Psi(\lambda\omega,\eps,y)\tilde\beta_\lambda(\omega,\eps,x,y),
$$
 are bounded by $M\delta$, where $M>0$ is independent of $\lambda,\omega,\eps,x,y.$
\end{remark}



\begin{lemma}\label{spectral-gap} Under assumption (\ref{A-stability}), there exists $\lambda_0>0$ such that 
$$
   \|A_\lambda^{-1}\|\le \mu^{-1},\quad 0<\lambda\le \lambda_0,\quad{\rm where}\ \mu=\dfrac{\|\beta_y(0)\|+1}{2}.
$$
\end{lemma}
 
\noindent The proof comes by observing that $A_0$ is an identity matrix and by deriving 
\begin{equation}\label{q<1}
  \|B\|=\|\beta_y(0)\|\le q<1
\end{equation}
from (\ref{A-stability}).
\vskip0.2cm

\noindent  
Since $\mu$ of Lemma~\ref{spectral-gap} verifies $
  \|B\|<\mu 
$, Lemma~\ref{spectral-gap} ensures the spectral gap property required by Hadamard-Perron theorem (as stated in Theorem~6.2.8 of \cite{Katok}) when applied to the map $F_\lambda.$ 

\vskip0.2cm

\noindent When combined together, conclusions of Corollary~\ref{Glambdaproof} and Lemma~\ref{spectral-gap} ensure that conditions of  Hadamard-Perron theorem (Theorem~6.2.8 of \cite{Katok}) hold and the following result about the existence of invariant curve for $G_\lambda$ takes place.

\begin{corollary}\label{corHP} Let assumptions 
of Theorem~\ref{thm1} hold and  
let $\lambda_0$ be the minimum of the two $\lambda_0$ given by Corollary~\ref{Glambdaproof} and Lemma~\ref{spectral-gap}. Then, for any $\lambda\in(0,\lambda_0],$ there exists a unique $\psi_\lambda\in C^1(\mathbb{R}\times\mathbb{R}\times\mathbb{R}^{k_1},\mathbb{R}^{k_2})$ such that 
\begin{equation}\label{conccorHP}
   G_\lambda({\rm graph}(\psi_\lambda))={\rm graph}(\psi_\lambda)\quad {and}\quad \|D\psi_\lambda\|\le 1.
\end{equation}
\end{corollary}



\noindent As the first two lines of $G_\lambda$ are identities, Corollary~\ref{corHP} implies the following.


\begin{corollary}\label{LemmaNew} In the settings of Corollary~\ref{corHP}, the conclusion~(\ref{conccorHP}) implies that 
$$
  G_\lambda(\omega,\eps,{\rm graph}(\psi_\lambda(\omega,\eps,\cdot)))=(\omega,\eps,{\rm graph}(\psi_\lambda(\omega,\eps,\cdot))),\ \ \omega,\eps\in\mathbb{R},\  \lambda\in(0,\lambda_0].
$$
\end{corollary}


\vskip0.2cm

\noindent 
In what follows, $\psi_\lambda$ will also refer to the one given by Corollary~\ref{corHP}.
We want to 
show that the graph of 
$
\psi_\lambda(\omega,\eps,\cdot)
$  is an invariant manifold also for $(x,y)\mapsto F_\lambda(\omega,\eps,x,y).$ Since $F_\lambda$ and $G_\lambda$ coincide on the domain where the bump function (\ref{bumpfunction}) equals 1, we now show that $(\omega,\eps,\psi_\lambda(\omega,\eps,x))$ belongs to this domain. 

\begin{lemma}\label{omega0}
    Let assumptions 
of Theorem~\ref{thm1} hold. Then there exists $\lambda_0>0$ such that for any $\lambda\in(0,\lambda_0]$ and for any $\omega\in\mathbb{R}$, if a $\varphi\in C^1(\mathbb{R}^{k_1},\mathbb{R}^{k_2})$ satisfies
$$
  G_\lambda(0,0,{\rm graph}(\varphi))=(0,0,{\rm graph}(\varphi)),
$$    
then $\varphi=0.$
\end{lemma}

\noindent {\bf Proof.} Let $(x,y)$ be such that $y=\varphi(x).$ Then
\begin{equation}\label{00}
  (0,0,x,\bar y)=G_\lambda(0,0,x,y),
\end{equation}
and by Lemma~\ref{LemmaNew}, $\bar y=\varphi(x)$ yielding \begin{equation}\label{yielding}
y=\bar y.
\end{equation} 
On the other hand, (\ref{00}) implies
\begin{equation}\label{tmp} 
\bar y=
\beta_y(0)y+\Psi(0,0,y)\tilde\beta_\lambda(0,0,x,y).  
\end{equation}


\noindent Applying the Taylor Theorem as in formula (\ref{betaestimate}),  
we get 
$$
  \Psi(0,0,y)\tilde\beta_\lambda(0,0,x,y)=\Psi(0,0,y)\int_0^1 \left(\beta_y(0,0,x,\lambda ty)y-\beta_y(0)y\right)dt.
$$
Thus, since the support of $y\mapsto\Psi(0,0,y)$ is bounded, there exists $L>0$ such that 
$
\left\|\Psi(0,0,y)\tilde\beta_\lambda(0,0,x,y)\right\|\le L\|\lambda y\|,
$
where $L>0$ doesn't depend on $x\in\mathbb{R}^{k_1}$ and $y\in\mathbb{R}^{k_2}$.
Therefore, we can choose $\lambda_0>0$ sufficiently small to have 
$$
  \|\Psi(0,0,y)\tilde\beta_\lambda(0,0,x,y)\|\le \dfrac{1-q}{2}\|y\|,\quad  x\in\mathbb{R}^{k_1},\ y\in\mathbb{R}^{k_2},\ \lambda\in(0,\lambda_0].
$$
Therefore, if such an $\lambda_0$ is taken from the beginning, then (\ref{tmp}) implies
\begin{equation}\label{implies}
   \|\bar y\|\le \dfrac{q+1}{2}\|y\|.
\end{equation}
Relations (\ref{yielding}) and (\ref{implies}) imply  $y=\bar y=0.$ The proof of the lemma is complete.\qed

\vskip0.2cm

\noindent Lemmas~\ref{LemmaNew} and \ref{omega0} imply that the  curve $\psi$ given by Corollary~\ref{corHP} satisfies \begin{equation}\label{zero}\psi_\lambda(0,0,x)\equiv 0,\quad\lambda\in(0,\lambda_0],
\end{equation}
if $\lambda_0$ is taken as the minimum of the two $\lambda_0$ given by Corollary~\ref{corHP} and Lemma~\ref{omega0}.

\begin{corollary}\label{corcor}    Let all assumptions 
of Theorem~\ref{thm1} hold. Let $\lambda_0$ be taken as the minimum of the two $\lambda_0$ given by Corollary~\ref{corHP} and Lemma~\ref{omega0}. Then 
\begin{equation}\label{psiomega0}
  \|\psi_\lambda(\omega,0,x)\|=0,\quad \omega\in\mathbb{R},\ x\in\mathbb{R}^{k_1},\ \lambda\in(0,\lambda_0].
\end{equation}    
\end{corollary}

\noindent {\bf Proof.} Fix an  $\omega\in\mathbb{R}.$ Let
$
  M=\sup_{x\in\mathbb{R}^{k_1}}\|\psi_\lambda(\omega,0,x)\|.
$
By (\ref{conccorHP}) and (\ref{zero}),
$
  M\le |\omega|.
$
Fix an arbitrary $\bar x\in\mathbb{R}^{k_1}$ and put $\bar y=\psi_\lambda(\omega,0,\bar x)$. By Corollary~\ref{LemmaNew}, there exist $x\in\mathbb{R}^{k_1}$ 
such that 
\begin{equation}\label{replacedby}
  (\omega,0,\bar x,\bar y)=G_\lambda(\omega,0,x,y),
\end{equation}
where $y=\psi_\lambda(\omega,0,x).$
Same arguments that were used to get (\ref{implies}) from (\ref{00}) can be repeated when $\omega=0$ is replaced by an arbitrary $\omega$ because $y\mapsto\Psi(\omega,0,y)$ is still bounded. In other words, (\ref{replacedby}) still implies (\ref{implies}), i.e.
\begin{equation}\label{new2.14}
\|\psi_\lambda(\omega,0,\bar x)\|=  \|\bar y\|\le \dfrac{q+1}{2}\|y\|.
\end{equation}
Since, $\bar x\in\mathbb{R}^{k_1}$ was chosen arbitrary,  the estimate (\ref{new2.14}) yields
$
   M\le\dfrac{q+1}{2}M,
$
i.e. $M=0$. The proof of the lemma is complete.\qed

\vskip0.2cm

\noindent Due to estimate (\ref{conccorHP}), the statement (\ref{psiomega0}) extends to
\begin{equation}\label{psiomega0extends}
\|\psi_\lambda(\omega,\eps,x)\|\le |\eps|,\quad \omega,\eps\in\mathbb{R},\ x\in\mathbb{R}^{k_1},\ \lambda\in(0,\lambda_0],
\end{equation}    
that allows us to conclude that, for suitable $\omega$ and $\eps$, $\psi_\lambda(\omega,\eps,\cdot)$ is also an invariant curve for $(x,y)\mapsto F_\lambda(\omega,\eps,x,y).$ Specifically, we have the following corollary of Corollary~\ref{corcor}.

\begin{corollary}\label{corcorcor}
In the settings of  Corollary~\ref{corcor}, 
\begin{eqnarray}
  & & F_\lambda(\omega,\eps,{\rm graph}(\psi_\lambda(\omega,\eps,\cdot)))=(\omega,\eps,{\rm graph}(\psi_\lambda(\omega,\eps,\cdot))), \label{graphinvarianceF1}\\ & &\omega\in\left[0,\dfrac{1}{\lambda}\right],\ \eps\in\left(-\dfrac{r_1}{2},\dfrac{r_1}{2}\right),\ \|x\|\in\mathbb{R}^{k_1},\ \lambda\in(0,\lambda_0],\label{graphinvarianceF2}
\end{eqnarray}  
\end{corollary}

\noindent The proof of Corollary~\ref{corcorcor} comes by concluding
\begin{equation}\label{psi-estimate}
\|\psi_\lambda(\omega,\eps,x)\|\le\dfrac{r_1}{2},\quad\omega\in\mathbb{R},\  \eps\in\left(-\dfrac{r_1}{2},\dfrac{r_1}{2}\right),\ \|x\|\in\mathbb{R}^{k_1},\ \lambda\in(0,\lambda_0],
\end{equation}
from (\ref{psiomega0extends}) and by observing that 
$
F_\lambda(\omega,\eps,x,y)= G_\lambda(\omega,\eps,x,y),
$
for all $\omega,\eps,x$ from the domain (\ref{graphinvarianceF2}) and for $y
\in\dfrac{r_1}{2}\mathbb{D}^{k_2}.$

\begin{lemma}\label{lemma-h} The equality  
\begin{equation}\label{if1}
       F_1(\omega,\eps,x,y)=(\omega,\eps,\bar x,\bar y),
\end{equation}
holds if and only if
\begin{equation}\label{if2}
F_\lambda(h_\lambda(\omega,\eps,x,y))=h_\lambda(\omega,\eps,\bar x,\bar y),
\end{equation}
where
$
   h_\lambda(\omega,\eps,x,y)=\left(\dfrac{\omega}{\lambda},\dfrac{\eps}{\lambda^2},x,\dfrac{y}{\lambda}\right).
$
\end{lemma}

\noindent {\bf Proof.} Indeed, if (\ref{if1}) holds, then
\begin{equation}\label{if3}
\begin{array}{rcl}
\renewcommand{\arraystretch}{2} 
F_\lambda(h_\lambda(\omega,\eps,x,y))&=&\left(\begin{array}{c}
\dfrac{\omega}{\lambda},
\dfrac{\eps}{\lambda^2},
x+\omega\alpha\left(\omega,\eps,x,y\right),
\dfrac{1}{\lambda}\beta(\omega,\eps,x,y)\end{array}
\right)=\\
&=&\left(\begin{array}{c}
\dfrac{\omega}{\lambda},
\dfrac{\eps}{\lambda^2},
\bar x,
\dfrac{1}{\lambda}\bar y\end{array}
\right)=h_\lambda(\omega,\eps,\bar x,\bar y).
\end{array}
\end{equation}
Conversely, if (\ref{if2}) holds, then we get (\ref{if1}) by considering the middle equality of (\ref{if3}) and by multiplying the first, second, and forth components of this equality by $\lambda,$ $\lambda^2$, and $\lambda$ respectively.
The proof of the lemma is complete.
\qed

\vskip0.2cm

\noindent Lemma~\ref{lemma-h} allows to deform $\lambda$ in Corollary~\ref{corcorcor} to $\lambda=1$ and get the following.

\begin{corollary}\label{cor4}
    In the settings of Corollary~\ref{corcor}, 
\begin{equation}\label{F1}
  F_1\left(\omega,\eps,{\rm graph}\left(\lambda\psi_\lambda\left(\dfrac{\omega}{\lambda},\dfrac{\eps}{\lambda^2},\cdot\right)\right)\right)=\left(\omega,\eps,{\rm graph}\left(\lambda\psi_\lambda\left(\dfrac{\omega}{\lambda},\dfrac{\eps}{\lambda^2},\cdot\right)\right)\right),     
\end{equation}
for all 
\begin{equation}\label{domain}\omega\in[0,1],\ \ \eps\in \left(-\dfrac{r_1\lambda^2}{2},\dfrac{r_1\lambda^2}{2}\right),\ \ \lambda\in(0,\lambda_0].
\end{equation}
\end{corollary}

\noindent {\bf Proof.} By Corollary~\ref{corcorcor},
$$  F_\lambda^i\left(\dfrac{\omega}{\lambda},\dfrac{\eps}{\lambda^2},x,\psi_\lambda \left(\dfrac{\omega}{\lambda},\dfrac{\eps}{\lambda^2},x\right)\right)=\left(\dfrac{\omega}{\lambda},\dfrac{\eps}{\lambda^2},x_i,\psi_\lambda \left(\dfrac{\omega}{\lambda},\dfrac{\eps}{\lambda^2},x_i\right)\right),\quad i\in\{-1,1\},
$$
for all $\omega$, $\eps$, and $\lambda$ satisfying (\ref{domain}). Here $x_1$ and $x_{-1}$ of course depend on $x$, $\omega,$ $\eps,$ $\lambda.$ Applying Lemma~\ref{lemma-h},
$$  F_1^i\left(h_\lambda^{-1}\left(\dfrac{\omega}{\lambda},\dfrac{\eps}{\lambda^2},x,\psi_\lambda \left(\dfrac{\omega}{\lambda},\dfrac{\eps}{\lambda^2},x\right)\right)\right)=h_\lambda^{-1}\left(\dfrac{\omega}{\lambda},\dfrac{\eps}{\lambda^2},x_i,\psi_\lambda \left(\dfrac{\omega}{\lambda},\dfrac{\eps}{\lambda^2},x_i\right)\right),\ i\in\{-1,1\},
$$
which is the same as
$$  F_1^i\left(\omega,\eps,x,\lambda\psi_\lambda \left(\dfrac{\omega}{\lambda},\dfrac{\eps}{\lambda^2},x\right)\right)=\left(\omega,\eps,x_i,\lambda\psi_\lambda \left(\dfrac{\omega}{\lambda},\dfrac{\eps}{\lambda^2},x_i\right)\right),\quad i\in\{-1,1\}.
$$
The property obtained coincides with (\ref{F1}). The proof is complete. \qed 

\vskip0.2cm

\noindent \begin{remark}\label{C1C3}  \rm  
Corollary~\ref{cor4} implies that $\varphi_{\omega,\eps}$ defined as
\begin{equation}\label{requiredpsi}
\varphi_{\omega,\eps}(x)=\lambda_0\psi_{\lambda_0}\left(\dfrac{\omega}{\lambda_0},\dfrac{\eps}{\lambda_0^2},x\right)
\end{equation}
verifies property (\ref{C1}) of Theorem~\ref{thm1} with \begin{equation}\label{eps0}\eps_0=\dfrac{r_1\lambda_0^2}{2},
\end{equation}
see Remark~\ref{Ff} for the relation between $F_1$ and $f_{\omega,\eps}$.
 Such a $\varphi_{\omega,\eps}$ also verifies (\ref{C3}) of Theorem~\ref{thm1} by Corollary~\ref{corcor} and  inequality (\ref{psiomega0extends}). \end{remark}
 
 \noindent The rest of this section discovers that $\varphi_{\omega,\eps}$ fulfills (\ref{C2}) and the uniqueness property, where  a suitable definition of $r_0$ will be needed. In what follows, $[F_\lambda(\omega,\eps,x,y)]_{3-4}$ stays for the $x$ and $y$ outputs of $F_\lambda$, i.e. if $(\omega,\eps,\bar x,\bar y)=F_\lambda(\omega,\eps,x,y)$, then $[F_\lambda(\omega,\eps,x,y)]_{3-4}=(\bar x,\bar y).$ By analogy, $[F_\lambda(\omega,\eps,x,y)]_3=\bar x.$
 

\begin{lemma}\label{Lemma2.5new} Let assumptions 
of Theorem~\ref{thm1} hold. Then, there exists $\lambda_0>0$ such that, for all $\lambda\in(0,\lambda_0]$, $\omega,\eps,x_0$ from the domain (\ref{graphinvarianceF2}), and $y_0
\in\dfrac{r_1}{2}\mathbb{D}^{k_2},$
 the trajectory $\{(x_j,y_j)\}_{j=0}^\infty$ of $(x,y)\mapsto [F_\lambda(\omega,\eps,x,y)]_{3-4}$ 
is defined for all $j\in\mathbb{N}$ and 
\begin{equation}\label{satisfies1}
   \|y_j-\psi_\lambda(\omega,\eps,x_j)\|\to 0\quad{\rm as}\quad j\to\infty.
\end{equation}   
Furthermore,  the convergence is uniform with respect to $(x_0,y_0).$  
\end{lemma}

\noindent {\bf Proof.} We first establish an auxiliary relation. For an arbitrary $\tilde x\in\mathbb{R}^{k_1}$, let $\tilde y=\psi_\lambda(\omega,\eps,\tilde x).$ Then, by the invariance of $\psi_\lambda$ (Corollary~\ref{corcorcor}), the point $(\omega,\eps,\bar x,\bar y)=F_\lambda(\omega,\eps,\tilde x,\tilde y)$ satisfies $\bar y=\psi_\lambda(\omega,\eps,\bar x).$ Therefore,
\begin{eqnarray*}
& & \dfrac{1}{\lambda}\beta(\lambda\omega,\lambda^2\eps,\tilde x,\lambda\tilde y)=\bar y=\psi_\lambda(\omega,\eps,\bar x)=\psi_\lambda\left(\omega,\eps,\tilde x+\lambda\omega\alpha(\lambda\omega,\lambda^2\eps,\tilde x,\lambda \tilde y)\right)=\\
&& 
=\psi_\lambda\left(\omega,\eps,\tilde x+\lambda\omega\alpha(\lambda\omega,\lambda^2\eps,\tilde x,\lambda\psi(\omega,\eps,\tilde x))\right).
\end{eqnarray*}
Using this auxiliary relation,
\begin{eqnarray*}
&&    \|y_j-\psi_\lambda(\omega,\eps,x_j)\|\le \left\|y_j-\dfrac{1}{\lambda}\beta\left(\lambda\omega,\lambda^2\eps,x_{j-1},\lambda\psi_\lambda(\omega,\eps,x_{j-1})\right)\right\|+\\
&&+\left\|\underline{\dfrac{1}{\lambda}\beta\left(\lambda\omega,\lambda^2\eps,x_{j-1},\lambda\psi_\lambda(\omega,\eps,x_{j-1})\right)}-\psi_\lambda(\omega,\eps,x_j)\right\|=\\
&&=\left\|\dfrac{1}{\lambda}\beta(\lambda\omega,\lambda^2\eps,x_{j-1},\lambda y_{j-1})-\dfrac{1}{\lambda}\beta\left(\lambda\omega,\lambda^2\eps,x_{j-1},\lambda\psi(\omega,\eps,x_{j-1})\right)\right\|+\\
&&+\left\|\underline{\psi_\lambda\left(\omega,\eps,x_{j-1}+\lambda\omega\alpha(\lambda\omega,\lambda^2\eps,x_{j-1},\lambda\psi_\lambda(\omega,\eps,x_{j-1}))\right)}-\right.\\
&&\left.-\psi_\lambda\left(\omega,\eps,x_{j-1}+\lambda\omega\alpha(\lambda\omega,\lambda^2\eps,x_{j-1},\lambda y_{j-1})\right)\right\|,
\end{eqnarray*}
where one pair of equal terms is underlined for reader's convenience.
By Lemma~\ref{magnitude}, there exists $\lambda_0>0$ such that
\begin{eqnarray*}
&& \hskip-0.6cm \left\|\dfrac{1}{\lambda}\beta(\lambda\omega,\lambda^2\eps,x,\lambda y_1)-\dfrac{1}{\lambda}\beta\left(\lambda\omega,\lambda^2\eps,x,\lambda y_2\right)\right\|=\\
&& \hskip-0.6cm=\left\| \tilde\beta_\lambda(\omega,\eps,x,y_1)-
\tilde\beta_\lambda(\omega,\eps,x, y_2) +\beta_y(0)y_1-\beta_y(0)y_2\right\|\le \left(q+\dfrac{1-q}{3}\right)\|y_1-y_2\|,
\end{eqnarray*}
for all $\omega\in\left[0,\dfrac{1}{\lambda}\right],$ $\eps\in\left(-\dfrac{r_1}{2},\dfrac{r_1}{2}\right),$ $x\in\mathbb{R}^{k_1}$, $\|y_1\|,\|y_2\|\le \dfrac{r_1}{2}.$
For this domain, we denote by $L>0$ the Lipschitz constant of $y\mapsto\alpha(\lambda\omega,\lambda^2\eps,x,y)$, which is bounded by (\ref{A-regularity}). Recalling that by (\ref{conccorHP}) the Lipschitz constant of $x\mapsto \psi_\lambda(\omega,\eps,x)$
doesn't exceed 1, we conclude
$$
\|y_j-\psi_\lambda(\omega,\eps,x_j)\|\le \left(q+\dfrac{1-q}{3}\right)\|y_{j-1}-\psi_\lambda(\omega,\eps,x_{j-1})\|+L\lambda\|\psi_\lambda(\omega,\eps,x_{j-1})-y_{j-1}\|.
$$
Reducing $\lambda_0>0$, if necessary, to have $L\lambda_0\le \dfrac{1-q}{3},$ we get
$$
\begin{array}{l}
  \|y_j-\psi_\lambda(\omega,\eps,x_j)\|\le \left(q+\dfrac{2}{3}(1-q)\right)\|y_{j-1}-\psi(\omega,\eps,x_{j-1})\|\le...\\
  ...\le \left(q+\dfrac{2}{3}(1-q)\right)^j\|y_0-\psi(\omega,\eps,x_0)\|,
\end{array}
$$
from where both $y_j\in r_1\mathbb{D}^{k_2}$ (see (\ref{psi-estimate})) and the convergence conclusion follow. \qed

\begin{corollary}\label{corconvergence} Let assumptions 
of Theorem~\ref{thm1} hold. Let $\lambda_0>0$ be given by Lemma~\ref{Lemma2.5new} and $\lambda\in(0,\lambda_0],$ $\omega\in[0,1],$ $\dfrac{\eps}{\lambda^2}\in\left(-\dfrac{r_1}{2},\dfrac{r_1}{2}\right).$ Then the trajectory $\{(x_j,y_j)\}_{j=0}^\infty$ of $(x,y)\mapsto [F_1(\omega,\eps,x,y)]_{3-4}$ with the initial condition $(x_0,y_0)\in\mathbb{R}^{k_1}\times \dfrac{\lambda r_1}{2}\mathbb{D}^{k_2}$ is defined, 
\vskip-0.35cm
$$
  \left\|y_j-\lambda\psi_\lambda\left(\dfrac{\omega}{\lambda},\dfrac{\eps}{\lambda^2},x_j\right)\right\|\to 0\quad{\rm as}\ \ j\to\infty,
$$
  and the convergence is uniform with respect to $(x_0,y_0).$ 
\end{corollary}

\noindent {\bf Proof.} 
By Lemma~\ref{Lemma2.5new}, the trajectory $\{\tilde x_j,\tilde y_j\}_{j=0}^\infty$ of $(x,y)\mapsto \left[F_\lambda\left(\dfrac{\omega}{\lambda},\dfrac{\eps}{\lambda^2},x,y\right)\right]_{3-4}$ with the initial condition $(\tilde x_0,\tilde y_0)\in \left(\mathbb{R}^{k_1},\dfrac{r_1}{2}\mathbb{D}^{k_2}\right)$ is defined and satisfies (\ref{satisfies1}). Therefore, the sequence $\{x_j,y_j\}_{j=0}^\infty$ defined by $(x_j,y_j)=(\tilde x_j,\lambda \tilde y_j)$ satisfies the given initial condition and   
$
   \left\|\dfrac{1}{\lambda}y_j-\psi_\lambda\left(\dfrac{\omega}{\lambda},\dfrac{\eps}{\lambda^2},x_j\right)\right\|\to 0,$  as $j\to 0.
$ By construction,
$F_\lambda\left(\dfrac{\omega}{\lambda},\dfrac{\eps}{\lambda^2},x_{j-1},\dfrac{y_{j-1}}{\lambda}\right)=\left(\dfrac{\omega}{\lambda},\dfrac{\eps}{\lambda^2},x_{j},\dfrac{y_{j}}{\lambda}\right).$ Therefore, $F_1(\omega,\eps,x_{j-1},y_{j-1})=(\omega,\eps,x_j,y_j)$ by Lemma~\ref{lemma-h}. The proof of the corollary is complete.          
\qed

\vskip0.2cm

\begin{remark}\label{remC2}
Corollary~\ref{corconvergence} concludes the validity of (\ref{C2}) of Theorem~\ref{thm1} for $\varphi_{\omega,\eps}$ given by (\ref{requiredpsi}) with $r_0$ defined as
\begin{equation}\label{r0}
  r_0=\lambda_0 r_1.
\end{equation}
\end{remark}

\noindent As (\ref{C1}) and (\ref{C3}) were established earlier  (see Remark~\ref{C1C3}), it remains to prove the uniqueness of $\varphi_{\omega,\eps}$, where the crucial role is played by the following lemma.

\begin{lemma}\label{lem-Matthew} Let assumptions  
of Theorem~\ref{thm1} hold. Let $\lambda_0>0$ be given by Lemma~\ref{Lemma2.5new}, $\omega\in[0,1]$, and $\dfrac{\eps}{\lambda_0^2}\in\left(-\dfrac{r_1}{2},\dfrac{r_1}{2}\right).$ Then, for any $x\in\mathbb{R}^{k_1}$ and any $n\in\mathbb{Z}$, there exists $x_0\in\mathbb{R}^{k_1}$ such that the trajectory $\{(x_j,y_j)\}_{j\in\mathbb{Z}}$ of $(x,y)\mapsto \left[F_1(\omega,\eps,x,y)\right]_{3-4}$ with the initial condition $(x_0,0)$ is defined and
\begin{equation}\label{xnx}
   x_n=x.
\end{equation}
\end{lemma}

\noindent {\bf Proof.} Let $F_1^n$ be the $n$-th iteration of $F_1$. It is sufficient to prove that
\begin{equation}\label{need}
\left[F_1^n(\omega,\eps,\mathbb{R}^{k_1},0)\right]_3=\mathbb{R}^{k_1}.
\end{equation}
We have
\begin{eqnarray*}
    \left[F_1(\omega,\eps,x,0)\right]_3&=&x+\omega\alpha(\omega,\eps,x,0),\\    \left[F_1^2(\omega,\eps,x,0)\right]_3&=&x+\omega\alpha(\omega,\eps,x,0)+\omega\alpha\left(\omega,\eps,x+\omega\alpha(\omega,\eps,x,0),\beta(\omega,\eps,x,0))\right).
\end{eqnarray*}
Considering $\left[F_1^3(\omega,\eps,x,0)\right]_3$, the formula will again begin with $"x+"$ and will contain a sum of three $\omega\alpha(...)$ with various arguments inside that continuously depend on $x$. To summarize, iterating $n$ times, we  get 
\begin{equation}\label{R(x)}\left[F_1^n(\omega,\eps,x,0)\right]_3=x+R(x),
\end{equation}
where $R$ is bounded on $\mathbb{R}^{k_1}$ by (\ref{A-regularity}). This representation implies, for any $\|x\|>0$,
\begin{equation}\label{1st}
\dfrac{\left<\left[F_1^n(\omega,\eps,x,0)\right]_3,x\right>}{\|x\|}=\|x\|+\left<R(x),\dfrac{x}{\|x\|}\right>\ge \inf\limits_{x\in\mathbb{R}^{k_1}\backslash\{0\}}\left<R(x),\dfrac{x}{\|x\|}\right>
\end{equation}
and 
\begin{equation}\label{2nd}\left\|\left[F_1^n(\omega,\eps,x,0)\right]_3\right\|\to\infty\quad {\rm as} \quad \|x\|\to\infty,
\end{equation}
where the right-hand-side of (\ref{1st}) is an element of $\mathbb{R}$.
The combination of (\ref{1st}) and (\ref{2nd}) yields  (\ref{need}) according to  \cite[Theorem~2.2.4]{Mawhin}.  The proof of the lemma is complete. \qed

\begin{remark} Note,  \cite[Theorem~2.2.4]{Mawhin} uses the Euclidean norm, while the norms in (\ref{1st})-(\ref{2nd}) are not necessary Euclidean. However,  (\ref{1st})-(\ref{2nd}) also hold with respect to the Euclidean norm because of the equivalence of norms in $\mathbb{R}^{k_1}.$
\end{remark}

\begin{remark}\label{rem-periodicity}\rm Using function $R(x)$ and formula (\ref{R(x)}), we can formulate the following extended version of statement (\ref{xnx}) 
\begin{equation}\label{extended}
x_n=
x_0+R(x_0)=x.
\end{equation}
Assume now that $R(x)$ in (\ref{R(x)}) is $T$-periodic with respect to $k$-th argument, i.e.
$R(z+Te^j)=z$, $z\in\mathbb{R}^{k_1}$, where $e^j$ denotes the vector of $\mathbb{R}^{k_1}$ with $1$ at $j$-th component and 0 at other components.
Let $\{(x_j,y_j)\}_{j\in\mathbb{Z}}$ be the trajectory of $(x,y)\mapsto f_{\omega,\eps}(x,y)$ with the initial condition $(x_0+Te^j,0)$. Then statement (\ref{extended}) takes the form 
\begin{equation}\label{xnx'}  x_n=x_0+Te^j+R(x_0+Te^j)=x_0+Te^j+R(x_0)=x+Te^j.
\end{equation}
\end{remark}

\noindent We will take advantage of Remark~\ref{rem-periodicity} in the next section when establishing periodicity in the invariant manifold $\varphi_{\omega,\eps}.$



\begin{corollary}\label{cor-uniqueness}
    Let all assumptions of Theorem~\ref{thm1} hold and let $\varphi_{\omega,\eps}$ be as defined in (\ref{requiredpsi}) with $\omega\in[0,1],$ $\eps\in(-\eps_0,\eps_0)$ and with $\eps_0$ and $\lambda_0$ given by (\ref{eps0}) and Lemma~\ref{Lemma2.5new}. Let $r_0$ be given by (\ref{r0}). If $\varphi\in C^1(\mathbb{R}^{k_1},r_0\mathbb{D}^{k_2})$ satisfies (\ref{C1})-(\ref{C2}), then
    $     \varphi(x)=\varphi_{\omega,\eps}(x),\quad x\in\mathbb{R}^{k_1}.
    $
\end{corollary}

\noindent {\bf Proof.} Assume that 
\begin{equation}\label{assume}
\varphi(x)\not=\varphi_{\omega,\eps}(x)
\end{equation}
for some $x\in\mathbb{R}^{k_1}.$ By the uniformness of the convergence in (\ref{C2}), there exists $n\in\mathbb{N}$ such that
$$
  \|y_n-\varphi_{\omega,\eps}(x_n)\|\le{\|\varphi(x)-\varphi_{\omega,\eps}(x)\|}/{3},\quad
  \|y_n-\varphi(x_n)\|\le{\|\varphi(x)-\varphi_{\omega,\eps}(x)\|}/{3},
$$
for any trajectory $\{(x_j,y_j)\}_{j=0}^\infty$ of $f_{\omega,\eps}$ with initial condition $(x_0,y_0)\in\mathbb{R}^{k_1}\times r_0\mathbb{D}^{k_2}.$

\vskip0.2cm

\noindent Now we apply Lemma~\ref{lem-Matthew} with $x$ and $n$ defined above and obtain $x_0$ such that the trajectory $\{(x_j,y_j)\}_{j=0}^\infty$ of $f_{\omega,\eps}$ with the initial condition $(x_0,0)$ satisfies (\ref{xnx}). Therefore, we have
$$
\|  \varphi(x)-\varphi_{\omega,\eps}(x)\|=\|y_n-\varphi_{\omega,\eps}(x_n)-(y_n-\varphi_{\omega,\eps}(x_n))\|\le{2} \|  \varphi(x)-\varphi_{\omega,\eps}(x)\|/3,
$$
contradicting (\ref{assume}). The proof of the corollary is complete. \qed

\vskip0.2cm 

\noindent As statements (\ref{C1}), (\ref{C2}), (\ref{C3}) were established earlier, see Remarks~\ref{C1C3} and \ref{remC2},  Corollary~\ref{cor-uniqueness} completes the proof of Theorem~\ref{thm1}.

\section{Periodicity in the associated invariant manifold when $\alpha$ and $\beta$ are periodic with respect to one of the $x$ variables}\label{Additional}


\noindent In the next theorem  the upper index $x^j$ denotes the $j$-th component of $x\in\mathbb{R}^{k_1}$, while  $x_j$ denotes the $j$-th element of the sequence $\{x_j\}_{j\in\mathbb{Z}}.$

\begin{theorem}\label{corperiodicity} If assumptions (\ref{A-zero}) and (\ref{A-stability}) of Theorem~\ref{thm1} stay, but assumption (\ref{A-regularity}) holds in the stronger sense 
\begin{itemize}
\item[\ref{A-regularity}')] in addition to (\ref{A-regularity}),  the functions $x^k\to\alpha(\omega,\eps,x^1,...,x^k,...,x^{k_1},y)$ and $x^k\to\beta(\omega,\eps,x^1,...,x^k,...,x^{k_1},y)$ are $T$-periodic for some $1\le k\le k_1$,  all $x^j$, $j\not=k,$ and all $y\in\mathbb{R}^{k_2},$
\end{itemize}
then, in addition to the boundedness statement of (\ref{C3}), it holds that \\$x^k\to \varphi_{\omega,\eps}(x^1,...,x^k,...,x^{k_1})$ is $T$-periodic, for all $x^j$, $j\not=k.$
\end{theorem}


\noindent{\bf Proof.} As in Remark~\ref{rem-periodicity} we will use $e_j$ to denote the vector of $\mathbb{R}^{k_1}$ whose only nonzero component is the component $j$ and it equals 1. Assume that the statement of the corollary doesn't hold. Then, there exists $x\in\mathbb{R}^{k_1}$ such that
\begin{equation}\label{xdef}
   \Delta=\left\|\varphi_{\omega,\eps}(x)-\varphi_{\omega,\eps}(x+Te_j)\right\|\not=0.
\end{equation}
According to statement (\ref{C2}) of Theorem~\ref{thm1}, there exists $n\in\mathbb{N}$ such that
\begin{equation}\label{ndef}
   \|y_n-\varphi(x_n)\|\le{\Delta}/{3},   
\end{equation}
for any trajectory $\{(x_j,y_j)\}_{j\in\mathbb{Z}}$  of $f_{\omega,\eps}$ satisfying $(x_0,y_0)\in \mathbb{R}^{k_1}\times \{0\}$. 

\vskip0.2cm

\noindent Let $x_0$ be as given by Lemma~\ref{lem-Matthew} applied with $x$ and $n$ from (\ref{xdef}) and (\ref{ndef}) respectively. Consider two trajectories 
$\{(\hat x_j,\hat y_j)\}_{j\in\mathbb{Z}}$ and $\{(\tilde x_j,\tilde y_j)\}_{j\in\mathbb{Z}}$ of $f_{\omega,\eps}$ with the initial conditions
$
   (\hat x_0,\hat y_0)=(x_0,0)$ and $(\tilde x_0,\tilde y_0)=(x_0+Te_j,0).
$
By Lemma~\ref{lem-Matthew} and Remark~\ref{rem-periodicity}, it holds that
$
  \hat x_n=x,$ and $ \tilde x_n=x+Te_j.
$
Furthermore, using condition~(\ref{A-regularity}'), 
$
   \hat y_n=\tilde y_n.
$
Therefore,
$$
\Delta=\left\|\varphi_{\omega,\eps}(x)-\varphi_{\omega,\eps}(x+Te_j)\right\|=\left\|\hat y_n-\varphi_{\omega,\eps}(\hat x_n)-\left(\tilde y_n-\varphi_{\omega,\eps}(\tilde x_n)\right)\right\|\le{\Delta}/{3}+{\Delta}/{3},
$$
which implies $\Delta=0$
 contradicting  (\ref{xdef}). The proof of the corollary is complete.\qed

\section{Application to the existence of a closed invariant curve and one-dimensionality of the dynamics in perturbed hybrid systems} \label{sec-application}

\noindent Let $S$ be given by
$
  S=\{x\in\mathbb{R}^n:H(x)=0\}.
$
 In this section we apply the results of Section~\ref{Additional} in order to investigate the response of a limit cycle of hybrid system 
\begin{eqnarray}
&& \dot x=X(x), \label{np1} \\
 && x(t)=\Delta(x(t^-)), \quad{\rm if}\ \ x(t^-)\in S,\label{np2} 
\end{eqnarray}
to periodic perturbations.
Therefore, along with system (\ref{np1})-(\ref{np2}), we introduce the perturbed system
\begin{eqnarray}
&&\dot x=X(x)+\eps g(t,x,\eps), \label{ps1} \\
&& x(t)=\Delta(x(t^-)), \quad{\rm if}\ \ x(t^-)\in S, \label{ps2}
\end{eqnarray}
where $\eps$ is a parameter. 

\vskip0.2cm

\noindent Assuming that the right-hand-sides of (\ref{np1}) and (\ref{ps1}) are continuous, we denote by $t\mapsto \phi_{t,\tau,\eps}(v)$ the solution $x(t)$ of (\ref{ps1}) with the initial condition $x(\tau)=v.$ 

\begin{definition} Let $x\in\mathbb{R}^{k_2+1}$. We say that positive-time-to-return map $T_\eps$ is defined at $(\tau,x),$ if 
$\phi_{I,\tau,\eps}(\Delta(x))\cap S\not=\emptyset,$ where $I=[\tau,c)$ is the maximal interval of definition of solution $t\mapsto\phi_{t,\tau,\eps}(\Delta(x))$ (here $c$ is finite or $c=\infty$, see \cite[Theorem~3.1]{Hartman}). When $T_\eps(\tau,x)$ is defined, we put 
$T_\eps(\tau,x)=\inf\{t\in(\tau,c):\phi_{t^-,\tau,\eps}(\Delta(x))\in S\}.$
\end{definition}

\noindent Let $D:\mathbb{R}^{k_2+1}\mapsto \mathbb{R}^{k_2}$ be some parameterization of $S$ in the sense that 
$$   
D(r_1\mathbb{D}^{k_2})\subset S,\qquad\mbox{for some }r_1>0.
$$

\begin{definition}\label{Poi} For each $\tau\in\mathbb{R}$ and $u\in r_1 \mathbb{D}^{k_2}$, such that $T_\eps(\tau,D(u))$ is defined, we introduce the Poincar\'e map $P_\eps$ of time-dependent system (\ref{ps1})-(\ref{ps2}) as
\begin{equation}\label{Pformula1}
  P_\eps(\tau,u)=\left(
  T_\eps(\tau,\Delta(D(u))),D^{-1}\left(\phi_{T_\eps(\tau,\Delta(D(u))),\tau,\eps}(\Delta(D(u)))
  \right) 
  \right).
\end{equation}
\end{definition}
\begin{definition}\label{Poinp} For each $u\in r_1 \mathbb{D}^{k_2}$, such that $T(D(u))$ is defined, the Poincar\'e map $P$ of reduced system (\ref{np1})-(\ref{np2}) is
$$
  P(u)=
  D^{-1}\left(\phi_{T(\Delta(D(u))),0,0}(\Delta(D(u)))
  \right).
$$
\end{definition}

\noindent To derive $P(u)$ from $P_\eps(\tau,u)$, one first observes that because (\ref{ps1}) is autonomous when $\eps=0$, \begin{equation}\label{taking1}
T_0(\tau,x)=\tau+T_0(0,x).
\end{equation}Then, since both $\tilde x(t)=\phi_{t,0,0}(x)$ and $\hat x(t)=\phi(t)_{t+\tau,\tau,0}(x)$ are solutions of (\ref{np1}) and $\tilde x(0)=\hat x(0)=x,$ it holds that $\tilde x(t)\equiv \hat x(t)$ yielding
\begin{equation}\label{taking2}
   \phi_{T_0(\tau,x),\tau,0}(x)=\phi_{T_0(0,x)+\tau,\tau,0}(x)=\phi_{T_0(0,x),0,0}(x).
\end{equation}
Therefore, the $\bar u$ in $(\bar\tau,\bar u)=P_0(\tau,u)$ is independent of $\tau$ and is the required $P(u)$. 

\vskip0.2cm

\noindent Based on (\ref{taking1})-(\ref{taking2}) and using the Implicit Function Theorem, if 
\begin{itemize}
\item[\ref{A-regularity}'')]  $X,\Delta\in C^2(\mathbb{R}^{k_2+1},\mathbb{R}^{k_2+1}),$ $g\in C^2(\mathbb{R}\times \mathbb{R}^{k_2+1}\times\mathbb{R},\mathbb{R}^{k_2+1}$),  $H\in C^2(\mathbb{R}^{k_2+1},\mathbb{R}),$ $D\in C^2(\mathbb{R}^{k_2},\mathbb{R}^{k_2+1})$, $g(t+T_g,x,\eps)\equiv g(t,x,\eps)$ for some $T_g>0,$ 

\item[\ref{A-zero}'')] $P(0)$ is defined,  $P(0)=0$, i.e. $D(0)$ is the initial condition of a cycle of (\ref{np1})-(\ref{np2}), 
and
$$
  \nabla H(P(0))X(P(0))\not=0, 
$$
\end{itemize}
then there exist $r_1>0$ and $\eps_0>0$ such that $T_\eps(\tau,D(u))$ is defined and of class $C^2$ for all $\tau\in [0,T_g]$, $u\in r_1\mathbb{D}^{k_2}$, and $\eps\in(-\eps_0,\eps_0).$

\begin{lemma}\label{Teps} If (\ref{A-regularity}'')-(\ref{A-zero}'') hold then $T_\eps$ on $\mathbb{R}\times D(r_1\mathbb{D}^{k_2})$ is defined recursively from the definition on   $[0,T_g]\times D(r_1\mathbb{D}^{k_2})$ according to the formula 
$$
T_\eps(\tau\pm T_g,x)=T_\eps(\tau,x)\pm T_g,\quad \tau\in[\pm jT_g,\pm (j+1)T_g],\ x\in D(r_1\mathbb{D}^{k_2}), \ j=0,1,2,...
$$
\end{lemma}

\noindent Here $[-j T_g,-(j+1)T_g]$ should be understood as $[-(j+1) T_g,-j T_g]$. We will write a proof for just $j=0$. An extension to $j=\pm1,\pm2,...$ can be done by analogy.

\vskip0.2cm

\noindent {\bf Proof.}  
 Observe that both  $\tilde x(t)=\phi_{t,\tau,\eps}(x)$ and $\hat x(t)=\phi_{t+T_g,\tau+T_g,\eps}(x)$ are solutions of (\ref{ps1}). Indeed, by $T_g$-periodicity of $g$,
\begin{equation}\label{periodic3}
      \hat x'(t)=X(\hat x(t))+\eps g(t+T_g,\hat x(t),\eps)=X(\hat x(t))+\eps g(t,\hat x(t),\eps).
\end{equation}
Noting that $\tilde x(\tau)=\hat x(\tau)=x$, we conclude $\tilde x(t)\equiv \hat x(t)$, or
\begin{equation}\label{perprop}\phi_{t,\tau,\eps}(x)=\phi_{t+T_g,\tau+T_g,\eps}(x).
\end{equation}

\vskip0.2cm

\noindent Pick $(\tau,x)\in [0,T_g]\times D(r_1\mathbb{D}^{k_2}).$ Since $T_\eps(\tau,x)$ is defined, 
$$
\phi_{T_\eps(\tau,x)^-,\tau,\eps}(x)\in S,\quad \phi_{t,\tau,\eps}(x)\not\in S,\quad t\in (\tau,T_\eps(\tau,x)).
$$
By (\ref{perprop}), the latter relation can be rewritten as
$$
\phi_{(T_\eps(\tau,x)+T_g)^-,\tau+T_g,\eps}(x)\in S,\quad \phi_{t+T_g,\tau+T_g,\eps}(x)\not\in S,\quad t\in (\tau,T_\eps(\tau,x)),
$$
or, making the change of the variables $s=t+T_g$, as
$$
\phi_{(T_\eps(\tau,x)+T_g)^-,\tau+T_g,\eps}(x)\in S,\quad \phi_{s,\tau+T_g,\eps}(x)\not\in S,\quad s\in (\tau+T_g,T_\eps(\tau,x)+T_g),
$$
which means that $T_\eps(\tau+T_g,x)$ is defined and given by the formula $T_\eps(\tau+T_g,x)=T_\eps(\tau,x)+T_g$. The proof of the lemma is complete.\qed

\begin{remark}\label{periodicity2}
     Lemma~\ref{Teps} implies that $T_g$-periodicity of $g(t,x,\eps)$ in $t$ ensures $T_g$-periodicity of
$$
   \alpha(\tau,u,\eps)=T_\eps(\tau,D(u))-\tau
$$
in $\tau$. Furthermore, observe that property (\ref{perprop}) (established in the proof of 
Lemma~\ref{Teps}) implies that 
$
\phi_{\tau+s,\tau,\eps}(x)$
is $T_g$-periodic in $\tau$ for each fixed $s$. Therefore, the function $$\beta(\tau,u,\eps)=D^{-1}\left(\phi_{\tau+\alpha(\omega,\tau,u,\eps),\tau,\eps}(D(u))\right)$$ is $T_g$-periodic in $\tau$.
\end{remark}

\noindent To apply Theorem~\ref{corperiodicity} we rewrite (\ref{Pformula1}) in the form
\begin{equation}\label{Pformula2}
  P_\eps(\tau,u)=\left(
  \tau+\alpha(\tau,D(u),\eps),\beta(\tau,u,\eps)
  \right).
\end{equation}

\begin{theorem}\label{thmappl} Assume that the reduced system (\ref{np1})-(\ref{np2}) satisfies conditions (\ref{A-regularity}'') and (\ref{A-zero}''). Assume additionally that
\begin{itemize}

\item[\ref{A-stability}'')] the absolute values of the eigenvalues of $P'(0)$ are strictly between 0 and 1. 
\end{itemize}

\noindent Then there exist $\eps_0>0$, $r_0\in (0,r_1]$, and a norm $\|\cdot\|$ in $\mathbb{R}^{k_2}$ such that for any $\eps\in(-\eps_0,\eps_0)$ there exists a unique $\varphi_{\eps}\in C^1(\mathbb{R},r_0\mathbb{D}^{k_2})$ satisfying
\begin{enumerate}
\item[\ref{C1}'')]
$
P_{\eps}(\mathrm{graph}(\varphi_{\eps}))=\mathrm{graph}(\varphi_{\eps});
$
\item[\ref{C2}'')]

the trajectory $\{(\tau_j,u_j)\}_{j=0}^\infty$ of $P_\eps$ with initial condition $(\tau_0,u_0)\in\mathbb{R}^{k_1}\times r_0\mathbb{D}^{k_2}$ is defined for all $j\in\mathbb{N}$ and 
$$
   \|u_j-\varphi_{\eps}(\tau_j)\|\to 0\quad{\rm as}\quad j\to\infty,
$$
 uniformly with respect to the initial condition $(\tau_0,u_0).$ 
\end{enumerate}

Furthermore
\begin{enumerate}
\item[\ref{C3}'')]
$\varphi_{\eps}$ is $T_g$-periodic, $\varphi_{\eps}$ depends continuously on $\eps\in(-\eps_0,\eps_0)$,  and $\varphi_{0}(x)\equiv 0.$ 
\end{enumerate}

\end{theorem}

\noindent {\bf Proof.} The regularity properties of (\ref{A-regularity}'') imply those of (\ref{A-regularity}) by differentiability of solutions on parameters and initial conditions (see e.g. Hartman \cite[Theorem~4.1]{Hartman}) and by differentiability of the implicit function (see e.g. Zorich \cite[\S8.5.4]{Zorich}). The periodicity property of 
(\ref{A-regularity}'') imply periodicity of $\alpha(\tau,u,\eps)$ and $\beta(\tau,u,\eps)$ (i.e. ensure that (\ref{A-regularity}') is satisfied) as noted in Remark~\ref{periodicity2}. Relations (\ref{taking1})-(\ref{taking2}) ensure  independence of $\alpha(\tau,u,0)$ and $\beta(\tau,u,0)$ of time and (\ref{A-zero}'') ensures $\beta(\tau,0,0)\equiv 0$, i.e. Assumption~(\ref{A-zero}) holds. Assumption (\ref{A-stability}'') implies that a norm $\|\cdot\|$ in $\mathbb{R}^{k_2}$ can be chosen in such a way that $u\mapsto\beta(\tau,u,0)$ contracts in $r_1\mathbb{D}^{k_2}$ for $\tau\in[0,T_g]$ with some $q\in[0,1)$, if $r_1>0$ is reduced sufficiently. Therefore, 
$\eps_0>0$ can be reduced so that $u\mapsto\beta(\tau,u,\eps)$ contracts in $r_1\mathbb{D}^{k_2}$ for $\tau\in[0,T_g]$ and $\eps\in[-\eps_0,\eps_0].$ The contraction statement extends to $\tau\in\mathbb{R}$ by $T_g$-periodicity of $\beta(\tau,u,\eps)$ in $\tau.$ To conclude the verification of Assumption (\ref{A-stability}), it remains to note that the lack of zero eigenvalues for $P'(0)$ implies invertibility of $u\mapsto \beta(\tau,u,0).$   Therefore, the conclusion follows by Theorem~\ref{corperiodicity} (applied with an arbitrary fixed value of $\omega$).\qed

\section{Conclusions} \noindent We proved that the dynamics of perturbed  autonomous systems with an attracting limit cycle is one-dimensional not only for smooth differential equations, but also for differential equations which are discontinuous  along a co-dimension 1 surface. The particular application of this paper is to hybrid differential equations, but the general result (Theorem~\ref{corperiodicity}) can be applied to any  Poincar\'e map of form (\ref{Pformula0}). 

\vskip0.2cm

\noindent The stability result that we conclude for hybrid systems (Theorem~\ref{thmappl}) implies that stable autonomous walking (including not falling) persists under periodic perturbations. Such a persistence doesn't necessary hold for arbitrary time-dependent perturbations \cite{Tedrake,Su} and  was established in the earlier literature via additional assumptions of a dwell-time type \cite{Gregg,Veer}. The availability of an extra parameter $\omega$ in the formulation of Theorem~\ref{thm1} can be explored by future researchers to obtain continuous dependence of the stable invariant manifold on e.g. the frequency of the periodic perturbation.



%


\section*{Acknowledgments} We thank NSF for financial support. 

\vskip0.5cm

\appendix

\noindent {\large \bf Appendix}

\begin{theorem}\label{diffthm}
    Let $G:\mathbb{R}^n\to\mathbb{R}^n$ be of the form $G(x)=Ax+g(x),$ where $A$ is an invertible $n\times n$ matrix. If $g\in C^1(\mathbb{R}^n,\mathbb{R}^n)$ and there exists $q>0$ such that
    $
       \|g'(x)\|\le q<\|A^{-1}\|^{-1},$ 
       $ x\in\mathbb{R}^n,
    $
    then $G$ is a surjective $C^1$-diffeomorphism from $\mathbb{R}^n$ to $\mathbb{R}^n.$
\end{theorem}



\noindent To prove surjectivity one can transform the problem of solving $G(x)=y$ for $x$ into the fixed point problem for the map $F(x)=A^{-1}(y-g(x)),$ for which the existence of a unique unique fixed point follows from the Banach Contraction Principle (see e.g. Zorich \cite[\S9.7]{Zorich2}). The $C^1$-regularity of the inverse comes by applying the Inverse Function Theorem (see e.g. Zorich \cite[\S8.6.1]{Zorich}) at each $y\in\mathbb{R}^n.$

\end{document}